# An Improved Dung Beetle Optimizer for Random Forest Optimization

Lianghao Tan*
W. P. Carey School of Business
Arizona State University
Tempe, USA
* Corresponding author; ltan22@asu.edu

Xiaoyi Liu
Ira A. Fulton Schools of Engineering
Arizona State University
Tempe, USA
xliu472@asu.edu

Dong Liu
Department of Computer Science
Yale University
New Haven, USA
pikeliu.mlsys@gmail.com

Shubing Liu
Department of Computer Science
University of North Carolina
at Chapel Hill
Chapel Hill, USA
sliu2@unc.edu

Weixi Wu
Tandon School of Engineering
New York University
New York, USA
ww2147@nyu.edu

Huangqi Jiang
College of Computing
Georgia Institute of Technology
Atlanta, USA
jianghuangqi@gmail.com

*Abstract*— **To improve the convergence speed and optimization accuracy of the Dung Beetle Optimizer (DBO), this paper proposes an improved algorithm based on circle mapping and longitudinal-horizontal crossover strategy (CICRDBO). First, the Circle method is used to map the initial population to increase diversity. Second, the longitudinal-horizontal crossover strategy is applied to enhance the global search ability by ensuring the position updates of the dung beetle. Simulations were conducted on 10 benchmark test functions, and the results demonstrate that the improved algorithm performs well in both convergence speed and optimization accuracy. The improved algorithm is further applied to the hyperparameter selection of the Random Forest classification algorithm for binary classification prediction in the retail industry. Various combination comparisons prove the practicality of the improved algorithm, followed by SHapley Additive exPlanations (SHAP) analysis.**

## I. Introduction

In the field of optimization problem solving, intelligent optimization algorithms have gradually become important tools for tackling complex problems [1-2] and exploration from multiple directions has been done to optimize model efficiency [3-5]. And the optimized algorithms have been used on various fields such as medical analysis, text processing, augmented reality, finance and economy [7-16]. Bio-inspired algorithms, as a typical approach, are widely applied in various fields such as engineering, economics, control, and logistics [17]. These algorithms simulate biological behaviors or evolutionary mechanisms found in nature to explore global optimal solutions, demonstrating advantages like strong adaptability and robustness. In recent years, many bio-inspired swarm intelligence algorithms have been proposed and successfully applied, such as Ant Colony Optimization (ACO), Particle Swarm Optimization (PSO), and Dung Beetle Optimizer (DBO) [18-20]. Among these, DBO is a relatively novel optimization algorithm, inspired by the foraging and dung-rolling behavior of dung beetles in nature. It simulates the navigation, search, and positioning abilities of dung beetles to solve complex optimization problems.

With its unique mechanism, DBO shows advantages in handling high-dimensional optimization problems. Compared to traditional optimization algorithms, DBO better balances global exploration and local exploitation, achieving a dynamic balance between exploration and exploitation during the solving process [20].

However, despite DBO's outstanding performance in many applications, it has some noticeable shortcomings. First, DBO's convergence speed is relatively slow, especially when handling large-scale problems, where the global search phase tends to consume significant computational resources. Second, DBO's local search ability is limited, making it prone to getting trapped in local optima, which affects the optimization results. Additionally, the randomness in the behavior of dung beetles may lead to insufficient algorithm robustness, resulting in different outcomes across multiple runs and a lack of consistency. These issues limit the wider application of DBO in practical scenarios [20].

This paper proposes an Improved Dung Beetle Optimizer (IDBO) based on circular population initialization, aiming to enhance convergence speed, local search capability, and global stability. Traditional DBO algorithms use random initialization, which, while generating diverse initial solutions, can lead to uneven population distribution and longer convergence times. This paper introduces a geometrically based circular population initialization method to ensure uniform distribution of the initial population in the solution space, thereby improving search efficiency. Additionally, a cross-validation mechanism is introduced, which screens and evaluates the population through multiple rounds, enhancing the robustness of the algorithm.

Cross-validation not only reduces the likelihood of getting trapped in local optima but also selects high-quality individuals early in the iteration process, thereby improving overall optimization performance.

In the experimental section, multiple standard benchmark test functions are used to validate the improved CICRDBO, and its performance is compared to the classical DBO. Experimental results show that CICRDBO outperforms the original algorithm and other comparison algorithms in terms of convergence speed, accuracy, and robustness. Furthermore, CICRDBO is applied to a real engineering problem, where the improved algorithm is used for hyperparameter selection in the Random Forest classification algorithm to predict a binary classification problem in the retail industry, further verifying its effectiveness in handling complex real-world problems.

The main contributions of this paper are as follows: (1) A multi-strategy improved DBO is proposed, which significantly improves the initial diversity of the population, accelerates convergence speed, enhances the algorithm's global search capability, and enables faster convergence to the optimal solution. (2) The improved algorithm is applied to hyperparameter selection in Random Forest classification, and simulation experiments are conducted on retail industry data, followed by SHAP analysis.

## II. IMPROVED DBO-- CICRDBO

### A. Circle Mapping

In the DBO, the initial positions of the population are usually generated through random sampling. However, random sampling often results in an uneven distribution of sample points across the search space, which affects the convergence efficiency and accuracy of the algorithm. Therefore, this paper adopts a Circle mapping method to generate a more uniform and reasonable sample point distribution, thus enhancing the diversity of the initial population.

Circle mapping is a typical example of chaotic mapping, and its mathematical form is quite simple. Please see equation (1):

$$x_{k+1} = mod\left(x_k + b - \frac{a}{2\pi}sin(2\pi x_k), 1\right) \quad (1)$$

where the parameters are a = 0.5 and b = 0.2.

### B. Longitudinal-Horizontal Crossover Strategy

The longitudinal-horizontal crossover strategy is crucial for enhancing the performance of the DBO, addressing its tendency to get stuck in local optima by promoting broader exploration of the solution space. Longitudinal crossover enhances global search ability, helping the population escape local optima and discover better search directions, while horizontal crossover fine-tunes local areas, improving precision and ensuring critical regions are not overlooked. This combination boosts diversity, accelerates convergence, and, compared to the traditional DBO, significantly improves optimization accuracy and stability. Ultimately, the improved algorithm more efficiently and accurately finds the global optimum, reducing computational costs and enhancing overall performance.

### C. Horizontal Crossover Operation

The horizontal crossover operation in the DBO is similar to the crossover operation in genetic algorithms, primarily performing crossover calculations on the same dimensions of different individuals. To address the issue of the Dung Beetle Optimizer's insufficient global search ability, this paper applies a horizontal crossover strategy to optimize the positions of individuals in the population. First, the population individuals are randomly arranged, and crossover operations are performed on the $d_{th}$ dimension to update the positions of the individuals. This method enhances information exchange between individuals, improving the algorithm's global search capability and enabling the Dung Beetle Optimizer to more effectively escape local optima and search for the global optimum. Please see equations (2) and (3)

$$MSx_{i,d}^t = r_1 \times x_{i,d}^t + (1 - r_1) \times x_{j,d}^t + c_1 \times \left(x_{i,d}^t \times x_{j,d}^t\right) \quad (2)$$

$$MSx_{j,d}^t = r_2 \times x_{j,d}^t + (1 - r_2) \times x_{i,d}^t + c_2 \times \left(x_{j,d}^t \times x_{i,d}^t\right) \quad (3)$$

where $MSx_{i,d}^t$ and $MSx_{j,d}^t$ represent individuals generated through horizontal crossover from the vigilant d by $x_{i,d}^t$ and $x_{j,d}^t$ $r_1$ and $r_2$ are random numbers within the range [0, 1], and $c_1$ and $c_2$ are random numbers within the range [-1, 1].

In the DBO, after horizontal crossover operations, individuals are more likely to generate offspring within their respective hypercube spaces and on the outer edges, thereby expanding the search space and enhancing global search capability. The solutions generated by horizontal crossover need to be compared with the parent individuals, and those with higher fitness are retained. This means that the number of offspring generated in the outer edge space decreases linearly as the distance between the parent individuals increases. Under this mechanism, the algorithm can gradually converge to the optimal solution while ensuring precision and maintaining high convergence efficiency. The horizontal crossover operation effectively balances exploration and exploitation, significantly improving the overall performance of the DBO.

### D. Longitudinal Crossover Operation

In the later stages of the DBO, it tends to get stuck in local optima. The primary reason for this is that some individuals in the population prematurely reach local optima in certain dimensions, leading to faster convergence of the entire population without fully exploring the global optimal solution. After analysis, it was found that the DBO lacks the necessary mutation mechanism, making it unable to effectively intervene in individuals trapped in local optima, thereby limiting the algorithm's ability to continue approaching the global optimum In the later stages of the DBO, it tends to fall into local optima. The main reason for this is that some individuals in the population prematurely reach local optima in certain dimensions, causing the entire population to converge more quickly and fail to fully explore the global optimal solution. Analysis reveals that the DBO lacks the necessary mutation mechanism to effectively intervene in individuals that are already trapped in local optima, thus limiting the algorithm's ability to continue approaching the global optimum. Therefore, after performing horizontal crossover operations, it is necessary to further apply longitudinal crossover to the newly generated individuals to enhance the algorithm's ability to escape local optima. Through longitudinal crossover, the algorithm can expand the search range of individuals, further optimizing the population's global search performance, and avoiding the trap of local optima.

The crossover operation is performed on all dimensions of the newly generated individuals, with a lower probability than horizontal crossover, similar to the mutation operation in genetic algorithms. Assuming that the new individual $x_{i,d}^t$ generates offspring through longitudinal crossover in dimensions $d_1$ and $d_2$, the calculation method is as equation (4):

$$MSx_{i,d}^t = r_1 \times x_{i,d_1}^t + (1-r) \times x_{i,d_2}^t \qquad (4)$$

where $MSx_{i,d}^t$ represents the offspring individuals generated from individual $x_{i,d}^t$ through longitudinal crossover in the $d_1$ and $d_2$ dimensions, where $r \in [0,1]$. Similar to horizontal crossover, the offspring generated by longitudinal crossover must compete with the parent individuals, and the individuals with higher fitness are retained. Through this preferential selection mechanism, the dung beetle individuals involved in the crossover do not lose valuable dimensional information; instead, the diversity of the population is improved, and the quality of the solution is continuously enhanced. After performing longitudinal crossover, individuals that were previously trapped in local optima can fully utilize useful information from various dimensions, giving them a chance to escape local optima.

During the iteration process, if an individual escapes a local optimum in one dimension through longitudinal crossover, the improvement will quickly propagate through the entire population via horizontal crossover, solidifying the new solution in that dimension. This provides more opportunities for other dimensions trapped in local optima to escape. The combination of horizontal and longitudinal crossover operations effectively enhances the algorithm's convergence efficiency and accuracy in escaping local optima.

## III. METHODOLOGY

### A. Experimental Platform

The experimental environment includes a Windows 11 operating system with an Intel(R) Core(TM) i9-14700HX CPU and 16.0GB memory.

### B. Experimental Design

To validate the optimization performance and application value of the improved algorithm, this paper designed 10 benchmark function experiments. The population size is set to 30, and the maximum number of iterations is 500. Each experiment is run 30 times, and the optimal value, average value, and standard deviation are recorded.

### C. Experimental Results

To visually compare the two algorithms, the convergence effect graphs for the Sphere, Rastrigin, Ackley, and Griewank benchmark functions are shown in Figure1. From these graphs, it can be seen that the CICRDBO algorithm has a faster convergence speed on the Sphere and Ackley benchmark functions, and better convergence results on the Rastrigin and Griewank benchmark functions.

### D. Random Forest Classification

The data used for this experiment is sourced from a Kaggle dataset (https://www.kaggle.com/code/prashant111/xgboost-k-fold-cv-feature-importance/input) in the retail industry, containing 9 fields: Channel, Region, Fresh, Milk, Grocery,

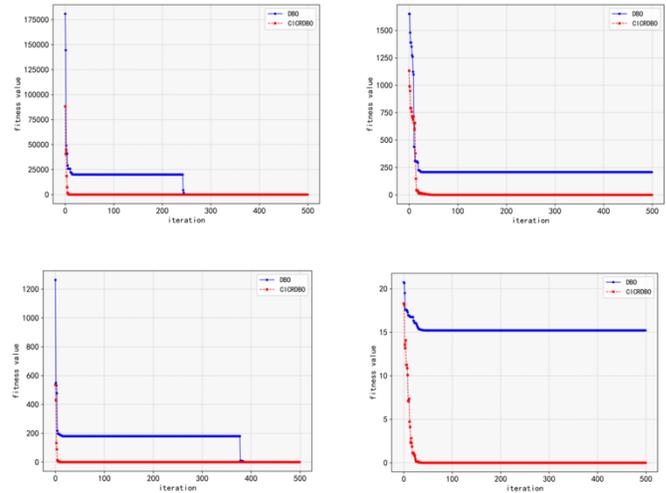

Figure 1. Experimental results on comparing algorithms

Frozen, Detergents Paper, and Delicatessen. Channel is the target variable, which predicts wholesale channels, while the other 8 are feature variables. Region is a discrete variable, while the rest are continuous variables.

First, statistical analysis was conducted on the data to observe their distribution. The proportions of Channel 1 and 2 are 67.7% and 32.3%, respectively, while the proportions of

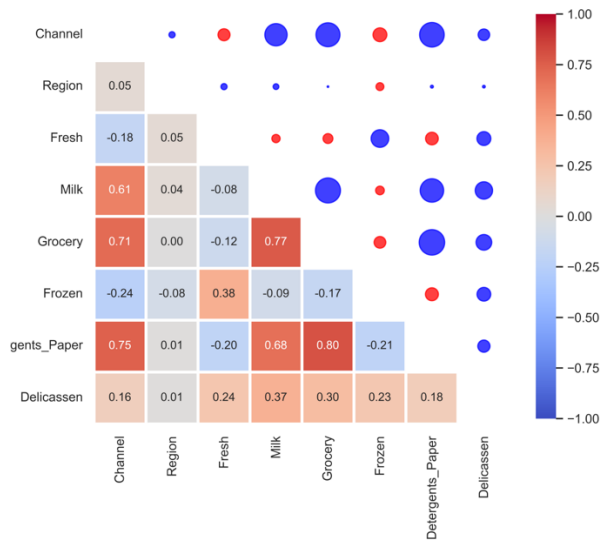

Figure 2. Correlations among feature variables and channels

Region 1, 2, and 3 are 17.5%, 10.7%, and 71.8%, respectively. Next, probability distribution for the continuous variables, shows that the six continuous variables generally follows a normal distribution. We then observed the relationship between feature variables and the target variable. When Region is 2, Channel 2 is more likely. Additionally, relationship between continuous variables and Channel indicates that when Channel is 1, the Fresh value is larger, while when Channel is 2, the Milk, Grocery, and Detergents Paper values are larger. A Spearman correlation was conducted to explore the relationship among feature variables and Channels. The variable with the highest correlation with Channel is Detergents Paper, with a correlation of 0.75, followed by Grocery with a correlation of 0.71. See Figure 2.

*E. Model Optimization*

A binary classification model was built using Random Forest. The dataset was split into 70% training and 30% testing. The model was trained on the training set and tested for generalization ability on the testing set using three methods: default hyperparameters, DBO-optimized hyperparameters, and CICRDBO-optimized hyperparameters. The results show that the CICRDBO-optimized hyperparameters achieved the best performance, as shown in Table 1.

TABLE I. MODELS FOR RANDOM FOREST CLASSIFICATION

| Models | Precision | Recall | F1 Score | AUC |
| --- | --- | --- | --- | --- |
| Default Parameters | 0.939 | 0.926 | 0.883 | 0.904 |
| DBO | **0.939** | 0.926 | 0.883 | 0.904 |
| CICRDBO | 0.946 | 0.950 | 0.883 | 0.915 |

*F. SHAP Analysis*

SHAP (SHapley Additive exPlanations) values are used to explain the results of machine learning models by calculating the contribution of each feature to the model's output [21]. The SHAP value is essentially a Shapley value from game theory, used to measure the contribution of each player in a game. In machine learning, SHAP values explain the importance of each feature for the prediction of each sample, helping us understand the model's decision-making process and uncover hidden patterns and trends in the data. SHAP values can be used to explain the prediction results of deep forest models. Results shows the feature importance based on average SHAP values, with the ranking of important features as follows: Detergents Paper (+0.24), Grocery (+0.15), Milk, Fresh (+0.05), Delicatessen (+0.01), Frozen (+0.01), and Region (+0.01).

A summary plot is a common visualization method that intuitively shows the contribution of each feature to the model's predictions. In Figure 3, each point represents a sample, the horizontal axis represents the SHAP value of the feature, and the vertical axis represents the feature. The color of the points indicates the value of the feature, with positive SHAP values indicating a positive influence and negative SHAP values indicating a negative influence.

While machine learning models are often considered black-box models that lack good interpretability, SHAP values effectively explain the prediction results. Figure 4 shows how a sample was predicted as Channel 1.

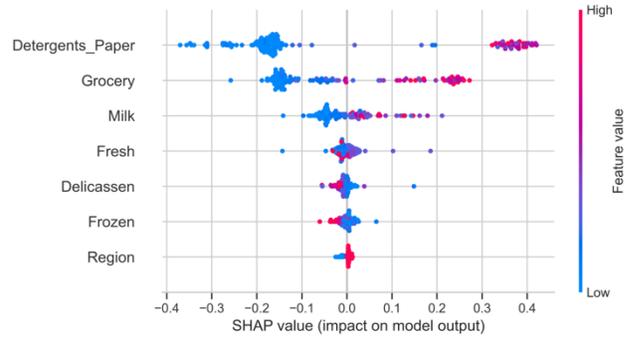

Figure 3. Summary plot of contribution of each feature

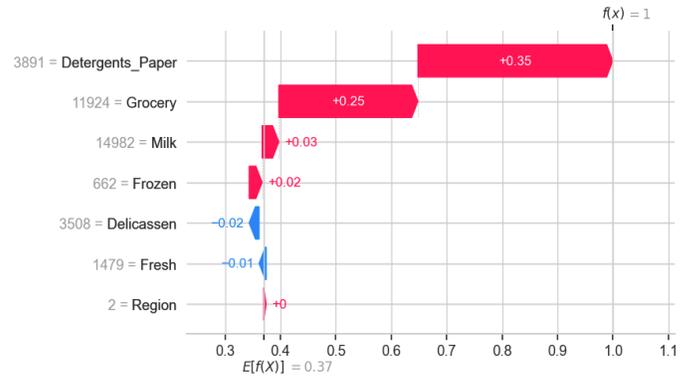

Figure 4. Example showing a sample predicted as Channel 1

*G. Sensitivity Analysis*

We conducted a sensitivity analysis to understand the effects of key CICRDBO parameters: Population Size, Max Iterations, Step Size, and Convergence Factor. Population Size: Affects diversity and computational load. Smaller sizes may limit global search, while larger sizes improve exploration but increase computation. Max Iterations: Fewer iterations may hinder full convergence; more iterations improve precision but require more computation. Step Size: Controls movement range. Larger steps risk skipping optimal solutions, while smaller steps may slow convergence. Convergence Factor: Balances global search and local refinement. A larger factor supports broader exploration, while a smaller factor enhances local accuracy but risks trapping in local optima.

Experiments show that adjusting these parameters within specific ranges significantly impacts performance. We recommend moderate population sizes and iteration counts, with step size and convergence factor tailored to balance exploration and precision, optimizing CICRDBO's performance for varied applications. This summary keeps the essential details while streamlining the content to meet brevity requirements.

## IV. Conclusion

In this paper, we proposed an improved DBO (i.e., CICRDBO) to address the limitations of the traditional DBO, specifically its slow convergence and tendency to get stuck in local optima. By incorporating circle mapping and a longitudinal-horizontal crossover strategy, the CICRDBO significantly enhances both the global search capability and convergence speed of the algorithm. Experimental results on multiple benchmark functions, including Sphere, Rastrigin, Ackley, and Griewank, demonstrate the superiority of the improved algorithm over the standard DBO in terms of optimization accuracy and efficiency.

Additionally, we applied the CICRDBO to optimize the hyperparameters of a Random Forest classifier, further validating its practical effectiveness in real-world applications. The results show that the improved algorithm not only achieves better performance in benchmark tests but also proves its utility in machine learning tasks. The CICRDBO has strong potential for further exploration and application in various optimization and classification problems, particularly those requiring efficient and accurate global search strategies.